\documentclass[12pt]{amsart}
\usepackage{amsmath}
\numberwithin{equation}{section}
\usepackage{amsthm}
\usepackage{amssymb}
\usepackage{amsfonts}
\usepackage{pinlabel}
\usepackage{graphicx}
\usepackage{relsize}
\newcommand{\nc}[2]{ \newcommand{#1}{#2} }

\nc{\avint}{ {- \hspace{-3.5mm} \int} }  % average integral

\nc{\R}{\mathrm{I \! R}}  % Real Numbers 
\nc{\N}{\mathrm{ I \! N}}  % Natural Numbers 

\newcommand{\tclosure}[1]{ \stackrel{\rule{.3 in}{.01 in}}{#1} }

\newcommand{\pclosure}[1]{ \stackrel{\rule{.5 in}{.01 in}}{#1} }

\DeclareRobustCommand{\rchi}{{\mathpalette\irchi\relax}}
\newcommand{\irchi}[2]{\raisebox{\depth}{$#1\chi$}}
\newcommand{\bigrchi}{\mathlarger{\mathlarger{\rchi}}}

\newcommand{\refeqn}[1]{ (\!\!~\ref{eq:#1}) } % gives references to
\newcommand{\refthm}[1]{ \!\!~\ref{#1} }    % equations or theorems

\nc{\Holder}{H\"{o}lder\ }

\nc{\ith}{ \ensuremath{\text{i}^{\text{th}}} }
\nc{\jth}{ \ensuremath{\text{j}^{\text{th}}} }
\nc{\kth}{ \ensuremath{\text{k}^{\text{th}}} }
\nc{\dst}{ \ensuremath{\text{1}^{\text{st}}_{\delta}} }
\nc{\dnd}{ \ensuremath{\text{2}^{\text{nd}}_{\delta}} }
\nc{\ost}{ \ensuremath{\text{1}^{\text{st}}} }
\nc{\tnd}{ \ensuremath{\text{2}^{\text{nd}}} }
\nc{\curl}{ \nabla \times }
\nc{\Div}{ \nabla \cdot }
\nc{\DC}{K}

\nc{\Ppl}{ \mathcal{M}^{+} }  \nc{\Pmn}{ \mathcal{M}^{-} }

\nc{\smiley}{ $\stackrel{\because}{\smile} \;$ }

\renewcommand{\eqref}[1]{(\ref{#1})}
\newcommand{\eat}[1]{}

\newcommand{\bsmall}{\begin{array}[c]{c}}
\newcommand{\esmall}{\end{array}}

\theoremstyle{plain}

\theoremstyle{definition}

\theoremstyle{plain}
\newtheorem{theorem}{Theorem}[section]

\newtheorem{corollary}[theorem]{Corollary}

\theoremstyle{definition}
\newtheorem{definition}[theorem]{Definition}
\newtheorem{remark}[theorem]{Remark}

\def\qed{\hfill\rule{1ex}{1ex}\\}
\newenvironment{pf}{\noindent {\bf Proof.}}{\qed}

\numberwithin{equation}{section}

\def\>{>_{\sigma}}

\title[Uniqueness of Mean Value Sets]{A Uniqueness Theorem for Mean Value Sets for Elliptic Divergence Form Operators}
\author[Armstrong]{Niles Armstrong}
\address{Milwaukee School of Engineering, Mathematics Department, 
1025 N. Broadway, 
Milwaukee, WI 53202}
 \email{armstrongn@msoe.edu}
\author[Blank]{Ivan Blank}
\address{Kansas State University, Mathematics Department, 
138 Cardwell Hall, 
Manhattan, KS 66506}
 \email{blanki@math.ksu.edu}

\begin{document}
\baselineskip 18pt

\begin{abstract}
We give background which shows the connection between the mean value theorem and the
obstacle problem, and then we prove that a set is a mean value set for an elliptic operator of the form
$Lu := \partial_i (a^{ij}(x) \partial_j u(x))$ if and only if it arises as the noncontact set of an obstacle
problem involving the Green's function of the operator.

\
 
 
\
 
\end{abstract}
\maketitle

\setcounter{section}{0}

\section{Introduction}   \label{Intro}

The Mean Value Theorem (MVT) for harmonic functions states that the value of a harmonic function $u$ at a point
$x_0$ can be recovered by taking an average of the function over any ball or sphere centered at that point.
(Appropriate inequalities can be shown for sub- and super-harmonic functions.)
This fact can then be used to develop the theory rather quickly; it provides quick and relatively simple proofs
of the strong maximum principle, the Harnack Inequality, and some of the fundamental a priori estimates which
lead to a variety of useful compactness properties.  (See \cite{GT, HL} for example.)  The
converse statement ``possessing the mean value property on every ball implies harmonicity'' is also known to
be true and is classical.

Historically, in spite of the importance of the the mean value theorem for the development of the theory for
the Laplacian, it has not been generalized to many operators.  On the other hand, in 1998 in the Fermi
Lectures on the Obstacle Problem, Caffarelli pointed out how the obstacle problem could be used to find
an analogous statement for a general uniformly elliptic divergence form operator of the form
$L := \partial_i a^{ij}(x) \partial_j$ and we give some of those details in Section \ref{MVTOP}.  He stated that
given an operator $L$ and a point
$x_0$ in the domain of $L,$ he could solve an obstacle problem in order to find an increasing family of
sets $\{ D_r(x_0) \}$ which were each comparable to $B_r(x_0)$ and so that if $v$ was any subsolution
to $Lv = 0,$ then the average
$$\int_{D_r(x_0)} \hspace{-.48in} \rule[.75mm]{.3cm}{.35mm} \ \ \ \ v(x) \; dx$$
would be increasing as a function of $r$ \cite{C2}.  Each of the $D_r(x_0)$ could be found as the unique noncontact
set for a solution to an obstacle problem.

Caffarelli stated the theorem more formally within a paper with Roquejoffre, but without giving the full proof,
he claimed that it could be found in the book on variational inequalities by Kinderlehrer and Stampacchia \cite{CR,KS}.
In fact, the theorem cannot be found there (or in other similar texts), so in \cite{BH1} the second author of this
paper and Z. Hao proved the theorem in detail.  Of course, even after stating what ``comparable to $B_r(x_0)$''
means precisely, it is still clear that the more that is known about the collection of these $D_r(x_0)$'s, the more
useful this theorem becomes, and since making this observation, both authors of the current paper have been
studying the properties of these sets \cite{Ar, ArmB, AryB, BBL}.    Of course, there is a natural definition to make:
\begin{definition}[Mean value sets]    \label{MVsets}
We will say that $D$ is a mean value set for the point $x_0$ and for the operator $L,$ if the mean value property
\begin{equation}
     u(x_0) = \int_{D} \hspace{-.21in} \rule[.75mm]{.3cm}{.35mm} \ \ \ \ u(x) \; dx %:= \frac{1}{|D|} \int_{D} u(x) \; dx
\label{eq:MVprop}
\end{equation}
holds for all $u$ satisfying $Lu = 0.$ (We will use the slash through the integral to denote division by the measure of
the domain of integration.) 
\end{definition}
\noindent
With this definition one can wonder if there are ever any mean value sets which do \text{not} arise as a noncontact
set of the obstacle problem that Caffarelli described.  The negative
answer is given in the main theorem of this paper which we state here:
\begin{theorem}[All mean value sets arise as noncontact sets from obstacle problems]   \label{AMVSAaNSfOP}
Assume that $D_{x_0} \! \subset \subset \R^n$ is an open connected mean value set for $x_0$ and for the operator $L,$ and
assume that $x_0 \in D_{x_0}.$  Then there is a choice of $r > 0$ so that 
\begin{equation}
\text{up to sets of measure zero} \ D_{x_0} \ \text{is equal to} \ D_r(x_0) \ \text{as given in Theorem\refthm{KeyMVT}\!\!.}
\label{eq:TheFuckingPoint}
\end{equation}
In other words, there exists a domain $\Omega$ with
$D_{x_0} \subset \Omega$ and a Green's function for $L$ on $\Omega$ and a value $r > 0$ so that
$D_{x_0}$ is the noncontact set (up to sets of measure zero again) obtained by solving the obstacle problem given
in Equation\refeqn{wrobprob}\!\!.
\end{theorem}
\noindent
In 1972, Kuran showed that if $D$ is an open connected mean value set for the point $x_0$ and for the Laplacian,
and if $x_0 \in D,$ then (except for sets of measure zero) $D$ must be an open ball centered at $x_0$ \cite{Ku}.
(See also some earlier related work by
Epstein, by Epstein and Schiffer, and by Goldstein and Ow \cite{Ep, EpSchiff, GO}.)  
Since the noncontact set for the specific obstacle problem that Caffarelli describes is always a ball when
$L = \Delta,$ the main new theorem in this paper can be viewed as an extension of the theorem due to Kuran.

There is also a connection with our work to certain results within the theory of quadrature domains, and we will be 
making use of some of the definitions and techniques found within that field.  We will
define quadrature domains below in Definition\refthm{ShapQDWS}\!\!, but for now we will say that mean value sets
are a specific type of quadrature domain and the solution to the obstacle problem is an analogue to the
modified Schwarz potential (see Definition\refthm{DefMSP}\!\!) which is associated to any quadrature domain.
We mention in particular a result due to Shahgholian, where he shows uniqueness of certain quadrature domains
under a variety of constraints.  Shahgholian does not consider quadrature domains for L-harmonic functions, but he
is dealing with much more general measures than simply multiples of the delta function.  (When viewed
from the perspective of quadrature domains, the measure associated to a mean value set is always a
multiple of the delta function.)  It is also worth observing that in order to prove his uniqueness theorem, he replaces
the assumption of nonnegativity of the modified Schwarz potential with a weaker condition that required a certain sum
to be nonnegative \cite{Shah}.  In our paper for the cases that we are considering, we are eventually able to show
that our analogue of the modified Schwarz potential is always nonnegative.

The proof of our main theorem will follow a part of Harold Shapiro's text on the Schwarz function somewhat closely,
and he in turn has followed Epstein and Schiffer fairly closely in the place where we are most dependent on his text
\cite{EpSchiff}.  Our ability to adapt this proof to L-harmonic functions is only possible because of certain recently
proven theorems on mean value sets which we state in Section \ref{MVTOP}.  The relevant material on Schwarz potentials
and quadrature domains can be found mostly
in \cite[Chapter 4]{Sh}, and near the end we borrow some terminology found within Sakai's text \cite{S}.  Having made
these citations, we note that Shapiro relies on fundamental solutions and convolutions, and we are forced to
rely on the Green's function and a slightly different construction in order to get to our analogue of the modified Schwarz
potential.  Furthermore, when Shapiro is studying quadrature domains, he is never studying $L$-harmonic
functions, but rather, he is considering harmonic functions, subharmonic functions, or complex analytic functions.  The
authors wish to thank Tim Mesikepp for pointing them in the direction of quadrature domains and Dave Auckly for a
helpful discussion.

\section{Connecting the Mean Value Theorem to the Obstacle Problem}    \label{MVTOP}

In this section we will explain how the MVT is related to the obstacle problem via the proof given by
Caffarelli in the Fermi lectures \cite{C2}.  The proof of the MVT for harmonic and/or subharmonic functions is
usually accomplished by computing:
\begin{equation}
       \frac{\partial}{\partial r} \left( 
               \int_{\partial B_r(x_0)} \hspace{-.54in} \rule[.75mm]{.3cm}{.35mm} \ \ \ \ \ \ \ \ u(x) \; dS_x
                                               \right)
\label{eq:UglyRoute}
\end{equation}
and then making use of a change of variables, the divergence theorem, and the assumptions on the function
$u.$  In the Fermi Lectures, however, Caffarelli gave a different proof based on producing a 
key test function to plug into the definition of ``weakly harmonic'' and/or ``weakly subharmonic'' function
which has a number of benefits over the one described above.  We recall the definition of
weakly harmonic:
\begin{definition}[Weakly harmonic and weakly subharmonic]   \label{WeaklyH}
Given a domain $\Omega \subset \R^n,$ we say that $u \in L^{1}_{loc}(\Omega)$ is \textit{weakly harmonic} if
for any $\phi \in C_c^{1,1}(\Omega),$ we have the following equality:
\begin{equation}
       \int_{\Omega} u \Delta \phi = 0 \;.
\label{eq:WeaklyH}
\end{equation}
We say that $u \in L^{1}_{loc}(\Omega)$ is \textit{weakly subharmonic} if
for any $\phi \in C_c^{1,1}(\Omega)$ with $\phi \geq 0$ we have the following inequality:
\begin{equation}
       \int_{\Omega} u \Delta \phi \geq 0 \;.
\label{eq:WeaklySubH}
\end{equation}
(Our ``$c$'' subscript means that our test functions are compactly supported within $\Omega.$)
\end{definition}
\noindent
%This definition arises naturally as the Euler-Lagrange equation that one finds when minimizing the Dirichlet
%integral (after one additional integration by parts), and of course, it allows us to discuss a notion of a solution to a second order
%equation that does not necessarily possess two classical derivatives.
%(Of course, in the case of weak variants
%of the Laplace Equation, Weyl's Lemma guarantees that weak solutions are always classical.)

Turning toward the relevant proof, we define $\Gamma(x)$ to be $|x|^{2-n}$ when the dimension $n$ is greater
than $2,$ and $-\ln |x|$ when $n = 2,$ and observe that because it is a multiple of the fundamental solution for the 
Laplacian, $\Delta \Gamma(x)$ vanishes outside of the origin.
Caffarelli creates a crucial auxiliary function by touching $\Gamma$ from below with a parabola of the form
$\alpha - \beta |x|^2.$  Because of the symmetries involved, it is clear that the ``touch'' occurs on a
sphere centered at the origin, and by varying the constants involved, he can make the touch happen on
a sphere with any desired radius.  His auxiliary function is defined to equal the parabola inside the sphere
where the touch occurs and equal to $\Gamma$ outside.  Following the proof he gives with only slight modifications,
we let $\psi_s(x)$ be the function created in this way when the touch occurs on
$\partial B_s.$  Now observe that the Laplacian of this function is a negative constant inside $B_s$ and it vanishes
outside of $B_s.$ In fact, because of the tangential touch, the function is $C^{1,1}$ and therefore it will not
``pick up a distribution'' along $\partial B_s.$  The upshot is that we can say that 
\begin{equation}
     \Delta \psi_s(x) = -C(s) \bigrchi_{B_s} \ \ \ \text{in} \ \ \R^n
\label{eq:LapOfHalf}
\end{equation}
in the weak sense.  Now the key test function is created as a difference of two of these auxiliary
functions.  In other words, let $0 < r < s$ and define
\begin{equation}
    \Phi_{r,s}(x) := \psi_r(x) - \psi_s(x) \;.
\label{eq:KeyTestFct}
\end{equation}
The function $\Phi_{r,s}$ is obviously $C^{1,1},$ and it is a simple exercise to show that it is nonnegative.  Since
it is equal to $\Gamma - \Gamma$ outside of $B_s,$ it obviously has compact support and so it satisfies everything
required of the test functions that we are allowed to plug into Equation\refeqn{WeaklySubH}\!\!.
Taking $u$ to be weakly subharmonic and plugging $\Phi_{r,s}$ into Equation\refeqn{WeaklySubH}we compute:
\begin{alignat*}{1}
     0 &\leq \int_{\Omega} u \Delta \Phi_{r,s} \\
        &= \int_{B_r} u \Delta \psi_r - \int_{B_s} u \Delta \psi_s \\
        &= - C(r) \int_{B_r} u + C(s) \int_{B_s} u \;.
\end{alignat*}
Now the only thing really left to do is verify that $C(r) = C|B_r|^{-1},$ and Caffarelli points out
that that can be accomplished very quickly by observing that $u \equiv 1$ is harmonic so that the inequality
above becomes an equality for that function, and after doing that the theorem is proved.

It is right after this proof, however, where Caffarelli makes a much bigger observation which is
contained in a short little remark \cite{C2}.  Within this remark, he
observes that although it \textit{appears} like the
construction of the auxiliary functions $\psi_s$ was dependent on the symmetries and smoothness of
Laplace's equation, it is actually \textit{not} the case.
%By viewing $\psi_s$ as a solution to an appropriate
%obstacle problem the symmetry and smoothness of the Laplacian become irrelevant.
In fact, after rewriting Equation\refeqn{LapOfHalf}as
\begin{equation}
     \Delta \psi_s(x) = -C(s) \bigrchi_{ \{ \psi_s < \Gamma \} } \ \ \ \text{in} \ \ \R^n \;,
\label{eq:LapOfHalfBetter}
\end{equation}
Caffarelli observed that what is really required of the auxiliary functions so that their difference can make a
good test function is that they each satisfy Equation\refeqn{LapOfHalfBetter}and that the difference of any two
of them has a sign and has compact support.  So, he can find an analogue for \textit{any} divergence form operator
of the form $L := \partial_i a^{ij}(x) \partial_j$ by solving an appropriate obstacle problem.  He stated that
by solving the problem:
\begin{equation}
    LW_r = \partial_i a^{ij}(x) \partial_j W_r = r^{-n} \bigrchi_{ \{ W_r > 0 \} } - \delta_{x_0}
\label{eq:CaffFormulation}
\end{equation}
where $\delta_{x_0}$ is the delta-function at the point $x_0,$ he could give his MVT for general divergence
form equations, and here we give a slightly improved version of his statement:
\begin{theorem} [Mean value theorem for divergence form elliptic PDE]  \label{KeyMVT}
Let $Lu := \partial_i (a^{ij}(x) \partial_j u(x))$ in the open connected domain $\Omega,$
assume that the $a^{ij}$ have ellipticity constants $\lambda$ and $\Lambda,$ and assume
that with this operator $L$ and this domain $\Omega$ there is a Green's function, $G.$ 
Next, for any $x_0 \in \Omega,$ define $D_r(x_0)$ to be the noncontact set for the solution
to the obstacle problem:
\begin{equation}
      \begin{array}{rll}
             L u &\!\!\!= -r^{-n} \bigrchi_{\left\lbrace u < G(\cdot, x_0) \right\rbrace} \ \ \ & \text{in} \ \Omega \\
              u &\!\!\!\le G(\cdot, x_0) \ \ \ & \text{in} \ \Omega \\
              u &\!\!\!= 0 \ \ \ \ &\text{on} \ \partial \Omega
      \end{array}
\label{eq:wrobprob}
\end{equation}
and assume that $r_0 = r_0(x_0)$ is the supremum of $r > 0$ such that $D_r(x_0) \subset \subset \Omega.$  Then
$r_0 > 0$ and for $r \in (0, r_0]$ the sets $\{ D_r(x_0) \}$ have the following properties:
\begin{enumerate}
  \item If $0 < r < s \le r_0,$ then $D_r(x_0) \subset D_s(x_0).$
  \item We have the inclusions: $B_{cr}(x_0) \subset D_r(x_0) \subset B_{Cr}(x_0),$
           with $c,$ $C$ depending only on $n,$ $\lambda,$
and $\Lambda.$
  \item For any $v$ satisfying $Lv \ge 0$ and $0 < r < s \le r_0$, we have
       \begin{equation}
              v(x_0) \le \frac{1}{|D_r(x_0)|}\int_{D_r(x_0)} v \le  \frac{1}{|D_s(x_0)|}\int_{D_s(x_0)} v \;,
       \label{eq:MVTres}
       \end{equation}
            and if $Lv \le 0,$ then the inequalities in Equation\refeqn{MVTres}are reversed, and of course, this fact
            leads to equalities when $Lv = 0.$
  \item Finally, we have $|D_r(x_0)| = r^n$ and $| \partial D_R(x_0) | = 0.$
\end{enumerate}
\end{theorem}

\begin{remark}[Some conventions]   \label{Odd}
Although we follow conventions of the text by Gilbarg and Trudinger in most respects \cite{GT}, and we use the 
paper by Littman, Stampacchia, and Weinberger on Green's functions \cite{LSW}, we have two conventions worthy of
note:
\begin{enumerate}
    \item We take the analyst's Laplacian and not the geometer's Laplacian.
    \item We take our ``Green's function'' $G$ so that $LG = -\delta_{x_0},$ so that $G$
             has a positive and not a negative singularity at $x = x_0.$
\end{enumerate}
\end{remark}
\begin{remark}[Existence and discussion of solutions to Equation\refeqn{wrobprob}\!\!]    \label{YeahTheyExist}
The existence of solutions to Equation\refeqn{wrobprob}follows from the theorems within \cite{BH1}.  In terms of
understanding this equation, one can describe the solution
(heuristically) as being a membrane that was inflated from below the Green's function, and the parameter
$r$ is related to the pressure of the inflation.  This description is only heuristic, but suggests correctly that
the problem can be reformulated as a convex minimization problem within the calculus of variations.
\end{remark}
\begin{remark}[Green's functions are preferable to a fundamental solution]  \label{GFbeatsFS}
Although Caffarelli and Blank-Hao use the fundamental solution and work on $\R^n$ in \cite{C2} and \cite{BH1},
especially within \cite{BBL} it becomes clear that the best way to find the $D_r(x_0)$'s is by solving the problem
given in Equation\refeqn{wrobprob}above where a Green's function is used.  One might worry about the
effect that changing $\Omega$ for an operator defined on all of $\R^n$ might change the $D_r(x_0)$'s, but
a key point observed within \cite{BBL} is that it has no effect as long as both domains are ``big enough.'' 
The point is that changing the ``outer'' set, $\Omega,$ leads to a new Green's function and a new solution
$u$ that \textit{have both been altered by the exact same L-harmonic function,} thereby preserving the
noncontact set.
\end{remark}
\begin{remark}[The height function]   \label{height}
Because it will come up later, we note that the ``height function'' $w(x) := G(x,x_0) - u(x)$ satisfies the following:
\begin{equation}
      \begin{array}{rll}
             Lw &\!\!\!= r^{-n} \bigrchi_{\left\lbrace w > 0 \right\rbrace} - \delta_{x_0} \ \ \ & \text{in} \ \Omega \\
              w &\!\!\!\ge 0 \ \ \ & \text{in} \ \Omega \\
              w &\!\!\!= 0 \ \ \ \ &\text{on} \ \partial \Omega \;,
      \end{array}
\label{eq:whappyobprob}
\end{equation}
and it follows from Theorem\refthm{KeyMVT}that $r^{-n} = |\{ w > 0 \}|^{-1}.$
\end{remark}

\section{Quadrature Domains and the Construction of the Modified Schwarz $L$-Potential}   \label{ACF}

%EDITED HERE
We take $\Omega$ to be an open connected bounded smooth set, and in particular, it has a Green's
function. We let $F(\Omega)$ denote the class of 
$L^1(\Omega)$ functions which are solutions of $Lu = 0$ in $\Omega.$
Following Shapiro, we make the following definitions:
\begin{definition}[Quadrature domains]  \label{ShapQDWS}
We say that $D$ is a \textit{Quadrature domain in the wide sense} (QDWS)
if there exists a distribution $\mu$ with compact support in $D$
for which the quadrature identity
\begin{equation}
     \int_{D} u \; dx = \langle \mu , u \rangle
\label{eq:QI}
\end{equation}
holds for all $u \in F(\Omega).$  We will say $D$ is a quadrature domain and
Equation\refeqn{QI}a quadrature identity if Equation\refeqn{QI}holds for a distribution $\mu$ whose
support is a finite set.
\end{definition}

Shapiro goes on to prove in \cite[Theorem 4.1]{Sh} that if $D$ is a QDWS, then it has a
\textit{modified Schwarz potential}.  From our point of view, these functions will correspond to the ``height
functions'' that can be produced as in Remark\refthm{height}\!\!.  Shapiro considers only the Laplacian and makes use
of some of the special properties of the fundamental solution in order to prove his result.  Our analogous
theorem for the L-harmonic case reads as follows:

\begin{theorem}[Mean value sets yield weak modified Schwarz $L$-potentials]  \label{MVSYWMSP}
Assume that $D_{x_0}$ is a mean value set for $x_0$ for the operator $L$
and which obeys $D_{x_0} \subset \subset \Omega.$  In other words we assume that for all $u$
which obey $Lu = 0$ in $\Omega$ we have
\begin{equation}
      \int_{D_{x_0}} \hspace{-.325in} \rule[.75mm]{.3cm}{.35mm} \ \ \ \ u(x) \; dx = u(x_0) \;.
\label{eq:MVset}
\end{equation}
Then there exists a unique function $w$ which satisfies:
\begin{equation}
     \begin{array}{rll}
           Lw \!\!\!&= \displaystyle{ \frac{1}{|D_{x_0}|} \bigrchi_{\left\lbrace D_{x_0} \right\rbrace} - \delta_{x_0} } \ \ \ \
                        &\text{in} \ \Omega \\
 \ \\
            w \!\!\!&\equiv 0 \ \ \ \ \ &\text{on all of} \ D_{x_0}^c \;.
     \end{array}
\label{eq:OPHF}
\end{equation}
\end{theorem}
\begin{definition}[Weak modified Schwarz $L$-potential]   \label{DefMSP}
We call a function $w$ which satisfies Equation\refeqn{OPHF}a \textit{weak modified Schwarz $L$-potential}
associated to the set $D_{x_0}$ and the operator $L.$  Based on the uniqueness given in the theorem, it is fair to
call it \textit{the} weak modified Schwarz $L$-potential in our situation.
\end{definition}
\begin{remark}[Nonnegativity?]   \label{NQuest}
Note that we do not yet assert that this function $w$ is nonnegative.
\end{remark}

\begin{pf}
We define $\mu$ to be the restriction of Lebesgue measure to $D_{x_0},$ so
$d\mu := \bigrchi_{\left\lbrace D_{x_0} \right\rbrace} dx,$ and we let $G(z,x)$ denote the Green's function
for $L$ on the set $\Omega$ as guaranteed to exist by \cite{LSW}.  (We recall our sign convention where our
Green's function is positive, and $LG = -\delta_{x_0}.$)  We set
$$W(z) := |D_{x_0}| G(z,x_0) - \int_{\Omega} G(z,x) \; d\mu(x) \;,$$
and by our definition of $\mu(x),$ this definition immediately gives:
\begin{equation}
     W(z) = |D_{x_0}| G(z,x_0) - \int_{D_{x_0}} G(z,x) \; dx \;.
\label{eq:HappyDefinition}
\end{equation}
It follows from the fact that $LG(y, x_0) \equiv 0$ in $D_{x_0}$ whenever $y \in D_{x_0}^c$ along with
the assumed mean value property of such functions with respect to $D_{x_0}$ that $W(z) \equiv 0$
on all of $D_{x_0}^c.$  Now by invoking \cite[Theorem 6.1]{LSW} we can conclude (with our sign convention) that
$$LW(z) = \bigrchi_{\left\lbrace D_{x_0} \right\rbrace} - |D_{x_0}| \delta_{x_0} \;.$$
Division by $|D_{x_0}|$ gives us the desired function, and based on the fact that it satisfies
Equation\refeqn{OPHF}it is automatically unique.
\end{pf}

\begin{remark}[Slightly weaker than what Shapiro requires]   \label{NotQuiteButGoodEnough}
Shapiro required that the modified Schwarz potential vanish \textit{along with its gradient} on the
boundary of the domain in question.  We do not make that requirement, but we will not need it either.
In the case where that happens, we will follow Shapiro and refer to it as a ``modified Schwarz $L$-potential''
or a ``\textit{true} modified Schwarz $L$-potential'' when we want to emphasize that its gradient is also vanishing
on the boundary of the relevant set.
\end{remark}

\section{The Uniqueness Theorem}   \label{TUT}

We turn to the proof of Theorem\refthm{AMVSAaNSfOP}\!\!.

\begin{pf}
We assume that $D_{x_0} \subset \subset \R^n$ is an open mean value set for the point $x_0$ and for
the operator $L.$  Using part (4) of Theorem\refthm{KeyMVT}we know that there exists an
$r > 0$ such that
\begin{equation}
    |D_r(x_0)| = |D_{x_0}| \;.
\label{eq:VolumeIsRight}
\end{equation}
By recalling Theorem\refthm{MVSYWMSP}we know that there is a function $W$ which satisfies
Equation\refeqn{OPHF}and so
\begin{equation}
     \begin{array}{rll}
           LW \!\!\!&= \displaystyle{ \frac{1}{|D_{x_0}|} \bigrchi_{\left\lbrace D_{x_0} \right\rbrace} - \delta_{x_0} } \ \ \ \
                        &\text{in} \ \Omega \\
 \ \\
            W \!\!\!&\equiv 0 \ \ \ \ \ &\text{on all of} \ D_{x_0}^c \;.
     \end{array}
\label{eq:WisMSP}
\end{equation}
Similarly, the function $W_0$ defined to be the height function from the obstacle problem that
produces $D_r(x_0),$ must satisfy:
\begin{equation}
     \begin{array}{rll}
           LW_0 \!\!\!&= \displaystyle{ \frac{1}{|D_{r}(x_0)|} \bigrchi_{\left\lbrace D_{r}(x_0) \right\rbrace} - \delta_{x_0} } \ \ \ \
                        &\text{in} \ \Omega \\
 \ \\
            W_0 \!\!\!&\equiv 0 \ \ \ \ \ &\text{on all of} \ D_{r}(x_0)^c \;.
     \end{array}
\label{eq:W0isMSP}
\end{equation}
(See Remark\refthm{height}\!\!.)
Now although we do not know whether $W$ is nonnegative everywhere, we absolutely do know that $W_0 \geq 0,$ 
%that $W_0$ vanishes together with its gradient on $\partial D_r(x_0),$ 
and that $W_0 > 0$ in all of $D_r(x_0).$
(In fact, that is how $D_r(x_0)$ is defined!)  We define
\begin{equation}
     \Upsilon := |D_{x_0}|(W - W_0) \ \ \text{so that} \ \
    L \Upsilon = \bigrchi_{\left\lbrace D_{x_0} \right\rbrace} - \bigrchi_{\left\lbrace D_{r}(x_0) \right\rbrace}
\label{eq:UpDef}
\end{equation}
and assume for the sake of contradiction that $D_{x_0} \ne D_{r}(x_0),$ and that they differ by more than a set of
measure zero.  It follows from De Giorgi's
Theorem that $\Upsilon$ is continuous, and since it vanishes outside of a compact set, it must attain its maximum and
minimum values.  We claim that the maximum value must be strictly positive.  Indeed, since $|D_r(x_0)| = |D_{x_0}|$
and since these sets are both open and differing on a set of positive measure, we can conclude that
$D_{x_0} \setminus D_{r}(x_0)$ has positive measure, and we can find a point $y$ of $\partial D_{x_0}$ which is
outside of the closure of $D_r(x_0)$
and so that $B_s(y) \cap D_{x_0}$ has positive measure for any $s > 0.$  Now if $s$ is sufficiently small so that
$B_s(y) \cap D_r(x_0) = \emptyset,$ then we
know that $L \Upsilon = \bigrchi_{\left\lbrace D_{x_0} \right\rbrace}$ in $B_s(y)$ and $\Upsilon(y) = 0 - 0 = 0.$
Since $L \Upsilon$ equals one on a subset of $B_s(y)$ with positive measure, we conclude that there exists a
$y_1 \in \partial B_s(y)$ such that
\begin{equation}
    \Upsilon(y_1) = \max_{\tclosure{B_s(y)}} \Upsilon(x) > 0 \;.
\label{eq:Yay}
\end{equation}
Thus, there is a point $z$ so that
\begin{equation}
    \Upsilon(z) = \max_{\R^n} \Upsilon(x) \geq \Upsilon(y_1) > 0 \;.
\label{eq:Yay2}
\end{equation}
Now since $W$ vanishes on the complement of $D_{x_0}$ and since $W_0 \geq 0$ we can conclude that
$z \in D_{x_0}.$  Since $D_{x_0}$ is open, there is a ball $B_q(z) \subset D_{x_0}$ and on this ball we must have:
\begin{equation}
     L \Upsilon(z) = 1 - \bigrchi_{\left\lbrace D_{r}(x_0) \right\rbrace} \geq 0 \;.
\label{eq:Yay3}
\end{equation}
By the strong maximum principle, we conclude that $\Upsilon \equiv \Upsilon(z)$ within $B_q(z),$ and now
we can repeat this process ensuring that we can never reach a boundary of the set
where $\Upsilon = \Upsilon(z) > 0.$ Since that set is compactly contained within $\R^n,$ it must have 
a boundary and that gives us a contradiction.
\end{pf}

\begin{corollary}[Existence of true modified Schwarz $L$-potentials]  \label{ETMSP}
Every mean value set has a modified Schwarz $L$-potential which vanishes together with its gradient on the
boundary of the set in question and furthermore, this modified Schwarz $L$-potential is always nonnegative.
\end{corollary}

%The following observation may not have been made before in the case of the Laplacian, but 
%\begin{corollary}
%If $D$ is a connected mean value set for the point $x_0$ and $x_0 \in D,$ then $D$ is almost everywhere
%equal to an open set.
%\end{corollary}

\begin{remark}[Extension to Riemannian manifolds]    \label{HappyToExtend}
The proofs and theorems given above extend easily to the mean value sets for the Laplace-Beltrami operator on
Riemannian manifolds as given by Benson, Blank, and LeCrone \cite{BBL}.
\end{remark}

%\begin{remark}[Other possible improvements]    \label{WeakenHypotheses}
%It certainly seems like there should be a way to weaken the hypotheses requiring ``connected'' and ``open''
%to something more measure theoretic.  Indeed, the proof given certainly shows that something like a
%connected fractal shape with positive measure which is not equal to an open set except for variation by a set
%of measure zero cannot be a mean value set.  This sort of improvement
%should work for the Laplacian as well as for any divergence form operator and so it would lead to
%\textit{proving} that every mean value set with positive measure is equal to an open set except for a set of
%measure zero.
%\end{remark}

It is worth observing that although mean value sets can obviously differ on sets of measure zero, the
modified Schwarz $L$-potential is unique (at least after taking the continuous representative), and so if
$W$ is the modified Schwarz $L$-potential for a given mean value set, then we suspect that the ``best''
representative of the mean value set would be the one defined by $\{ W > 0 \}.$  Another candidate
for the ``best'' representative was given in Sakai's text by the ``areal maximal domain'' \cite{S}.
After we change from the complex numbers to $\R^n$ his definition becomes the following: When given a domain
$D \subset \R^n$ we define the ``areal maximal domain'' by:
\begin{equation}
    [D] := \{ \; y \in \R^n \; : \; |B_r(y) \setminus D| = 0 \ \text{for some} \ r > 0 \; \} \;.
\label{eq:arealmaximaldef}
\end{equation}
For comparison, and with thoughts of singular points and lower dimensional singular sets on our minds, we
observe that a set of the form $\{ W > 0 \}$ could have ``slits'' and points missing, whereas by applying
the definition in Equation\refeqn{arealmaximaldef}we would ``fill in'' those points.  (Think of a punctured or
``slit'' disk versus a disk without those subsets missing.)  Taking the areal maximal domain seems to be
throwing away those singular points, and there is currently lots of interest in exactly those points.  (See \cite{FS}
for very recent work devoted to an understanding of such points.)

To further understand our preference, it is worth examining some of the special types of points that 
Sakai singled out as being interesting in the tenth section of his text \cite[Section 10]{S}.   Within that section, he
considered ``increasing families of domains'' which he denoted by ``$\{ \Omega(t) \}.$''  He
set
\begin{equation}
    \text{disc}\{ \Omega(t) \} := \left[ \bigcup_{t > 0} \Omega(t) \setminus \Omega(0) \right]
                  \setminus \bigcup_{t > 0} \partial \Omega(t) \;,
\label{eq:discset}
\end{equation}
and he made the following definitions for points within
$$S_{\Omega} := \bigcup_{t > 0} \Omega(t) \setminus \Omega(0) \;.$$
\begin{definition}[Stationary and stagnant points]   \label{StaPts}
A \textit{stationary point} $p \in S_{\Omega},$ is a point which belongs to $\partial \Omega(t)$ for all $t$ in an open interval.  
A \textit{stagnant point} $p \in S_{\Omega},$ is a point where there exists a $t(p) \ge 0$ and an $\epsilon(p) > 0$ such that
$|B_r(p) \setminus \Omega(t)| > 0$ for all $r > 0$ and all $t < t(p),$ and such that
$|B_{\epsilon(p)}(p) \setminus \Omega(t)| = 0$ for every $t \ge t(p).$
\end{definition}

Combining the vocabular of Sakai that we have defined above, together with the results found within
\cite{ArmB}, \cite{AryB}, and \cite{BBL}, and with the agreement that we take the mean value sets given
as the positivity set of the height function, i.e. we take $\Omega(t) := D_t(x_0)$ and not $\Omega(t) := [D_t(x_0)],$
we have the following corollary:
\begin{corollary}[Foliation results]   \label{FolRes}
The set $\text{disc}\{ \Omega(t) \}$ is empty, and under the added assumption that $a^{ij} \in C^{1,1}$
(or assuming that we are dealing with the Laplace-Beltrami operator on a Riemannian manifold)
the union of the connected components of the sets $\{ \partial D_t(x_0) \}$ give all of the leaves of a singular
foliation of $\pclosure{D_{r_0}(x_0)} \setminus \{ x_0 \}.$  In particular, in this case, the set of stationary points is
empty.  (We say ``singular foliation'' because some of the boundary components may be nonsmooth or even lower
dimensional sets.)
\end{corollary}
Now if, instead, we work with the areal maximal versions of our mean value sets, i.e. we let
$\Omega(t) := [D_t(x_0)],$ then we cannot foliate our domain with the new $\partial \Omega(t)$
in general as we will throw away some of our singular points.  As an example of a situation where we could
simultaneously produce a stagnant point and a point that failed to belong to the union of the boundaries of the
areal maximal versions of our mean value sets, we can consider the situation on a Riemannian manifold with a
long ``tendril.''  By choosing an $x_0$ near the base of the tendril and an $t > 0$ that is not \textit{too} small,
we can produce a mean value set that wraps around the tendril.  Then by increasing $t$ until the moment when
$D_t(x_0)$ swallows up the tendril, we can produce the desired point.

%%%%%%%%%%%%%%%%%%%%%%%%%%%%%%%%%%%%%%%%%%%%%%

%%%%%%%%%%%%%%%%%%%%%%%%%%%%%%%%%%%%%%%%%%%%%%

\bibliographystyle{plain}
\bibliography{BIB}
\end{document}